\documentclass{amsart}

\usepackage{amsfonts}
\usepackage{amsmath}
\usepackage{amsthm}
\usepackage{amssymb}
\usepackage[english]{babel}
\usepackage{graphics}
\usepackage{latexsym}
\usepackage{longtable} \usepackage{mathrsfs}
\usepackage{graphicx}
\usepackage{wrapfig} 
\usepackage[pdftex,colorlinks=true]{hyperref}

\newtheorem{theorem}{Theorem}
\newtheorem{corollary}[theorem]{Corollary}
\newtheorem{lemma}[theorem]{Lemma}

\newtheorem{claim}[theorem]{Claim}
\newtheorem{example}[theorem]{Example}

\theoremstyle{definition}
\newtheorem{definition}[theorem]{Definition}
\newtheorem{remark}[theorem]{Remark}

\renewcommand{\S}{\mathcal{S}}

\newcommand{\A}{\textrm{A}}

\newcommand{\R}{\mathbb{R}}

\newcommand{\mS}{\mathbb{S}}

\newcommand{\noi}{\noindent}
\newcommand{\ms}{\medskip}

\newcommand{\al}{\alpha}
\newcommand{\be}{\beta}
\newcommand{\ga}{\gamma}
\newcommand{\de}{\delta}
\newcommand{\De}{\Delta}

\newcommand{\la}{\lambda}
\newcommand{\ka}{\kappa}
\newcommand{\Om}{\Omega}

\newcommand{\lharpoonup}{-\!\!\!\!\rightharpoonup}

\newcommand{\larrow}{\longrightarrow}
\newcommand{\ot}{\otimes}

\newcommand{\ri}{\rightarrow}

\newcommand{\p}{\partial}
\newcommand{\sub}{\subseteq}

\newcommand{\by}{\times}

\newcommand{\sgn}{\textrm{sgn}}
\newcommand{\ess}{\textrm{ess}}

\newcommand{\cof}{\textrm{cof}}

\newcommand{\bt}{\begin{theorem}}\newcommand{\et}{\end{theorem}}
\newcommand{\bd}{\begin{definition}}\newcommand{\ed}{\end{definition}}
\newcommand{\bl}{\begin{lemma}}\newcommand{\el}{\end{lemma}}
\newcommand{\beq}{\begin{equation}}\newcommand{\eeq}{\end{equation}}
\newcommand{\bc}{\begin{claim}}\newcommand{\ec}{\end{claim}}
\newcommand{\bex}{\begin{example}}\newcommand{\eex}{\end{example}}
\newcommand{\bcor}{\begin{corollary}}\newcommand{\ecor}{\end{corollary}}
\newcommand{\bp}{\begin{proof}}\newcommand{\ep}{\end{proof}}

\newcommand{\BPL}{\medskip \noindent \textbf{Proof of Lemma} }

\newcommand{\BPT}{\medskip \noindent \textbf{Proof of Theorem} }

\numberwithin{equation}{section}

\begin{document}

\title[Existence and Uniqueness for Fully Nonlinear Systems]{Existence and Uniqueness of Global Solutions to Fully Nonlinear First Order Elliptic Systems}

\author{Nikos Katzourakis}
\address{Department of Mathematics and Statistics, University of Reading, Whiteknights, PO Box 220, Reading RG6 6AX, Berkshire, UK}
\email{n.katzourakis@reading.ac.uk}


\date{}



\begin{abstract} Let $F : \mathbb{R}^n \times \mathbb{R}^{N\times n} \rightarrow \mathbb{R}^N$ be a Carath\'eodory map. In this paper we consider the problem of existence and uniqueness of weakly differentiable global strong a.e.\ solutions $u: \mathbb{R}^n \longrightarrow \mathbb{R}^N$ to the fully nonlinear PDE system
\[\label{1} \tag{1}
F(\cdot,Du ) \,=\,  f, \ \ \text{ a.e.\ on }\R^n,
\]
when $ f\in L^2(\R^n)^N$. By introducing an appropriate notion of ellipticity, we prove existence of solution to \eqref{1} in a tailored Sobolev ``energy" space (known also as the J.L.\ Lions space) and a uniqueness a priori estimate. The proof is based on the solvability of the linearised problem by Fourier transform methods and a ``perturbation device" which allows to use of Campanato's notion of near operators, an idea developed for the 2nd order case.
\end{abstract}

\maketitle

\section{Introduction} \label{section1}

Let $n,N\geq 2$ and let also
\[
F\ :\ \R^{n} \by \R^{N \by n} \larrow  \R^N,
\]
be a Carath\'eodory map, namely
\[
\left\{
\begin{array}{l} 
x\mapsto F(x,Q) \text{ is measurable, for every } Q\in \R^{N\by n},\ms\\
Q\mapsto F(x,Q) \text{ is continuous, for almost every } x\in \R^{n}.
\end{array}
\right.
\]
In this paper we consider the problem of existence and uniqueness of global strong a.e.\ solutions $u : \R^n \larrow \R^N$ to the fully nonlinear PDE system
\beq  \label{1.1}
F(\cdot,Du ) \,=\, f, \ \ \text{ a.e.\ on }\R^n.
\eeq
To the best of our knowledge, the above problem has not been considered before in this generality. We will assume that our right hand side $f$ is an $L^2$ vector function, i.e.\ $ f\in L^2(\R^n)^N$. By introducing an appropriate ellipticity assumption on $F$, we will prove unique solvability of \eqref{1.1} for a weakly differentiable map $u$, together with a strong a priori estimate. In the above, $Du(x) \in \R^{N\by n}$ denotes the gradient matrix of of $u=u_\al e^\al$, namely $Du=(D_iu_a) e^\al \ot e^i$ and $D_i=\p/\p x_i$. Here and in the sequel we employ the summation convention when $i,j,k,...$ run in $\{1,...,n\}$ and $\al,\be,\ga,...$ run in $\{1,...,N\}$. Evidently, $\{e^i\}$, $\{e^\al\}$ and $\{e^\al \ot e^i\}$ denote the standard bases of $\R^n$, $\R^N$ and $\R^{N\by n}$ respectively. 

The simplest case of \eqref{1.1} is when $F$ is independent of $x$ and linear in $P$, that is when
\[
F(x,P)\, =\, \A_{\al \be j} P_{\be j}e^\al,
\]
for a linear operator $\A : \R^{N\by n}\larrow \R^N$. Then \eqref{1.1} becomes
\[
\A_{\al \be j} D_ju_\be \,=\, f_\al,
\]
which we will write compactly as
\beq \label{1.2}
\A: Du \,=\, f.
\eeq
The appropriate notion of ellipticity in this case is that the nullspace of the operator $\A$ contains no (non-trivial) rank-one lines. This means
\beq\label{1.3}
|\A :\eta \ot a | \, >\, 0,\ \ \ \eta \neq 0,  \ a \neq 0.
\eeq 
The prototypical example is the operator $\A : \R^{2 \by 2} \larrow \R^2$ given by
\[
\A \, =\, 
\left[
\begin{array}{cc|cc}
1 & 0 & 0 & 1\\
0 & \!\!\!-1 & 1 & 0 
\end{array}
\right]
\]
which corresponds to the Cauchy-Riemann differential operator. Linear elliptic first order systems with constant coefficients have been extensively studied in several contexts, since they play an important role in  Complex and Harmonic Analysis (see e.g.\ Buchanan-Gilbert \cite{BG}, Begehr-Wen \cite{BW}), Compensated Compactness and Differential Inclusions (Di Perna \cite{DP}, M\"uller \cite{Mu}), regularity theory of PDE (see chapter 7 in \cite{Mo} for Morrey's exposition of the Agmon-Douglis-Nirenberg theory) and Geometric Analysis with differential forms (Csat\'o-Dacorogna-Kneuss \cite{CDK}). 

The fully nonlinear case of \eqref{1.1} is much less studied. When $F$ is \emph{coercive instead of elliptic}, the problem is better understood. By using the analytic Baire category method of the Dacorogna-Marcellini \cite{DM} which is the ``geometric counterpart" of Gromov's Convex Integration, one can prove that, under certain structural and compatibility assumptions, the Dirichlet problem
\beq \label{1.4}
\left\{
\begin{array}{r}
F(\cdot,Du)\,=\, f, \ \ \ \ \text{ in }\Om, \smallskip\\
u\,=\, g,\ \, \text{ on }\p \Om,
\end{array}
\right.
\eeq
has \emph{infinitely many} strong a.e.\ solutions in the Lipschitz space, for $\Om \sub \R^n$ and $g$ Lipschitz. However, ellipticity and coercivity of $F$ are, roughly speaking, mutually exclusive. In particular, the Dirichlet problem \eqref{1.4} is not well posed when $F$ is either linear or elliptic. For example, the equation $u'-1=0$ has no Lipschitz solution on $(0,1)$ for which $u(0)=u(1)=0$. 

Herein we focus on the general system \eqref{1.1} and we consider the problem of finding an ellipticity condition which guarantees existence as well as \emph{uniqueness} of a strong a.e.\ weakly differentiable solution. In order to avoid the compatibility difficulties which arise in the case of the Dirichlet problem on bounded domains, we will consider the case of global solutions on the whole space. We will also restrict attention to the case of $f$ in $L^2(\R^n)^N$ and $n\geq 3$. This is mostly for technical simplicity and since the case of $n=2$ has been studied much more extensively. We will prove unique solvability of \eqref{1.1} in the ``energy" Sobolev space\footnote{We would like to thank the referee of this paper who pointed out to us that \eqref{1.5} is called in the literature the ``J.L.\ Lions space".}
\beq \label{1.5}
W^{1;2^*\!,2}(\R^n)^N\, :=\, \Big\{ u\in L^{2^*}(\R^n)^N\ \big| \ Du\in L^{2}(\R^n)^{Nn}\Big\} 
\eeq
where $n\geq 3$ and $2^*$ is the conjugate Sobolev exponent:
\[
2^*\, =\, \frac{2n}{n-2}.
\]
The technique we follow for \eqref{1.1} is based on the solvability of the linear constant coefficient equation \eqref{1.2} via the Fourier transform. In Section \ref{section2} we prove existence and uniqueness of a strong global solution $u\in W^{1;2^*\!,2}(\R^n)^N$ to \eqref{1.2}, for which we also have an \emph{explicit} integral representation formula for the solution (Theorem \ref{th1}). The essential idea in order to go from the linear to the fully nonlinear case is a perturbation device inspired from the work of Campanato \cite{C0}-\cite{C5} on the second order case of 
\beq \label{eq6}
F(\cdot,D^2u)\, = \, f.
\eeq
Campanato introduced a notion of strict ellipticity which requires that the nonlinear operator $F[u]:=F(\cdot,D^2u)$ is ``near" the Laplacian $\De u$ (see also Tarsia \cite{Ta1}-\cite{Ta3}). This implies unique solvability of \eqref{eq6} in $H^2\cap H^1_0$, by the unique solvability of the Poisson equation $\De u =f$ in $H^2\cap H^1_0$ and a fixed point argument. This notion can be also seen equivalently as a sort of strict pseudo-monotonicity, related to the Cordes condition (see Cordes \cite{Co1,Co2} and also Talenti \cite{T} and Landis \cite{L}). That such a stringent notion is required in order to guarantee well posedness is evident by counter-examples which are valid even in the linear scalar case of the second order elliptic equation 
\[
A_{ij}(x)D^2_{ij}u(x)\, =\, f(x) 
\]
with $A$ measurable and $A\geq \la I$ for $\la>0$ (see e.g.\ Pucci \cite{P}, and also Giaquinta-Martinazzi \cite{GM}). We also note that Campanato's notion has been weakened by Buica-Domokos in \cite{BD} to a notion of ``weak nearness", which still retains most of the features of (strong) nearness. In the same paper, the authors also use an idea similar to ours, namely a fully nonlinear operator being ``near" a general linear operator, but they implement this idea in the scalar case of 2nd order elliptic equations.

In Section \ref{section3} we introduce a notion of strict ellipticity which is inspired by Campanato's ellipticity, the latter being referred to as ``Condition A" in the literature (Definition \ref{def1}). Loosely speaking, our notion requires that $F$ is ``not too far" form a constant coefficient operator $\A$. We also introduce a related notion which we call pseudo-monotonicity and examine their connection (Lemma \ref{pr1}). Finally, we prove existence and uniqueness of solutions to \eqref{1.1} when $F$ is elliptic (Theorem \ref{th2}). This is based on the solvability of the linear system, our ellipticity assumption and the fixed point theorem in the guises of  Campanato's result of ``near operators" taken from \cite{C0}, which we recall herein for the convenience of the reader (Theorem \ref{th3}). A byproduct of our method is a strong uniqueness estimate in the form of a \emph{comparison  principle} for the distance of any solutions in terms of the distance of the right hand sides of the equations (Corollary \ref{cor1}).

\section{Existence-uniqueness-representation in the linear case} \label{section2}

In this section we prove existence and uniqueness of solutions to the linear constant coefficient system
\[
\A : Du \,=\, f, \ \ \text{ a.e.\ on }\R^n,
\]
when $\A : \R^{N \by n} \larrow \R^N$ is elliptic, that is the nullspace does not contain lines spanned by rank-one directions. This means that all rank-one lines are transversal to the nullspace:
\[
\A:\eta \ot a \,\neq \,0, \ \ \text{ when }\,  \eta\neq 0,\ a\neq 0.
\]
By compactness of the torus $\mS^{N-1}\by \mS^{n-1} \sub \R^N \by \R^n$, this is equivalent to
\beq \label{2.1}
|\A : \eta \ot a | \,\geq \, \nu\, |\eta| |a|, \ \ \eta \in \R^N, \, a\in \R^n,
\eeq
for some $\nu>0$, which can be taken to be the \emph{ellipticity constant of $\A$}:
\beq \label{2.1a}
\nu(\A)\ :=\ \min_{|\eta|=|a|=1} \big|\A : \eta \ot a \big|. 
\eeq
It is easy to see that \eqref{2.1} is equivalent to
\beq \label{2.2}
\min_{a \in \mS^{n-1}}\big| \det(\A a) \big| \, >\, 0,
\eeq
where $\A a$ is the $N \by N$ matrix
\[
\A a \, :=\, (\A_{\al \be j}\, a_j )\, e^\al \ot e^\be.
\]

\begin{example}
An example of $\A : \R^{2 \by 2} \larrow \R^2$ satisfying the above type of ellipticity is
\[
\A \, =\, 
\left[
\begin{array}{cc|cc}
\ka & 0 & 0 & \la\\
0 & \!\!\!-\mu & \nu & 0 
\end{array}
\right],
\]
where $\ka,\la,\mu, \nu>0$. A higher dimensional example of elliptic $\A : \R^{4\by 3}\larrow \R^4$ is given by
\[
\A \, =\, 
\left[
\begin{array}{rrr|rrr|rrr|rrr}
1 & 0 & 0    &      0 & \!\!\!-1 & 0    &    0 & 0 & \!\!\!-1   & 0& 0 & 0  \\
0 & 1 & 0    &      1 & 0 & 0    &     0 & 0 & 0   &  0& 0 & \!\!\!-1  \\
0 & 0 & 1    &      0 & 0 & 0    &    1 & 0 & 0   &   0& 1 & 0  \\
0 & 0 & 0    &      0 & 0 & 1    &    0 & \!\!\!-1 &0   &   1& 0 & 0  \\
\end{array}
\right]
\]
and corresponds to the electron equation of Dirac in the case where is no external force, the mass is zero and we are in the static case without time dependence. For a pair of complex functions $\phi,\psi : \R^3 \larrow \mathbb{C}$, the relevant system reads
\[
\left\{
\begin{split}
(D_1\, -\, i D_2)\phi\, +\, D_3\psi\, =\, F_1,\\
(D_1\, +\, i D_2)\psi\, -\, D_3\phi\, =\, F_2\\
\end{split}
\right.
\]
and in real coordinates becomes
\[
\left\{
\begin{split}
D_1u_1 + D_2u_2+D_3u_3 \phantom{+D_3u_3} & \ =f_1,\\
-D_2u_1 +D_1u_2 \phantom{+D_3u_3}+D_3u_4&\ =f_2,\\
-D_3u_1 \phantom{+D_3u_3}+D_1u_3-D_2u_4 &\ =f_3,\\
 \phantom{+D_3u_3}-D_3u_2+D_2u_3+D_1u_4&\ =f_4.
\end{split}
\right.
\]
In order to see that the tensor $\A$ is indeed is elliptic, we have that for any $a=(a_1,a_2,a_3)^\top\in \R^3$ with $a_1^2+a_2^2+a_3^2=1$, a quick calculation gives
\[
\A a \, =\, 
\left[
\begin{array}{rrrr}
a_1 & a_2 & a_3  &  0 \\
-a_2 & a_1 & 0  & a_3 \\
-a_3 & 0 & a_1  & -a_2 \\
0 & -a_3 & a_2  & a_1 \\
\end{array}
\right]
\]
and the rows (or the columns) of the $4\by4$ matrix $\A a$ form an orthomornal frame of vectors in $\R^4$. Hence, $\A a $ is an orthogonal matrix in $O(4,\R)$, from which it follows that $\det(\A a)=\pm 1$. 
\end{example}

The next theorem is the main result of this section. Before the statement and for the reader's convenience, we note that the elementary ideas of Fourier Analysis we use herein can be found e.g.\ in Folland \cite{F} and we follow the same notations as therein. In particular, for the Fourier transform and its inverse we use the conventions
\[
\widehat{u}(z)\,= \int_{\R^n}u(x)e^{-2\pi i x\cdot z}dx\ , \ \ \ \overset{\vee}{u}(x)\,= \int_{\R^n}u(z)e^{2\pi i x\cdot z}dz.
\]
Here ``$\cdot$" is the inner product of $\R^n$. With ``sgn" we denote the sign function on $\R^n$, that is $\sgn(x)=x/|x|$ when $x\neq 0$ and $\sgn(0)=0$. With ``$\cof(X)$" we denote the cofactor matrix of $X \in \R^{N \by N}$ and we will tacitly use the identity
\[
X\cof(X)^\top \ =\ \cof(X)^\top X \ = \ \det(X)I. 
\]

\bt[Existence-Uniqueness-Representation] \label{th1} Let $n\geq 3$, $N\geq 2$ and $\A : \R^{N \by n}\larrow \R^N$ a constant map satisfying \eqref{2.1}. Let also $f\in L^2(\R^n)^N$. Then, the problem
\[
\emph{\A} : Du \,=\, f, \ \ \text{ a.e.\ on }\R^n,
\]
has a unique solution $u$ in the space $W^{1;2^*\!,2}(\R^n)^N$ (see \eqref{1.5}), which also satisfies the estimate
\beq \label{2.3}
\|u\|_{L^{2^*}(\R^n)}\, + \ \|Du\|_{L^{2}(\R^n)}\, \leq\, C \|f\|_{L^{2}(\R^n)}
\eeq
for some $C>0$ depending only on $\emph{\A}$. 

Moreover, we have the following representation formula for the solution:
\beq \label{2.4}
u\, =\, -\frac{1}{2\pi i} \lim_{m\ri \infty}
\left\{  \widehat{h_m} \ast
\left[ 
\frac{\  \emph{\cof} \,(\emph{\A} \sgn )^\top}{ \det (\emph{\A} \sgn) \ } \overset{\vee}{f}
\right]^{\wedge} 
\right\}.
\eeq
In \eqref{2.4} $(h_m)^\infty_1 \sub \S(\R^n)$ is any sequence of even functions in the Schwartz class for which 
\[
\text{$0\, \leq\, h_m(x) \,\leq\, \frac{1}{|x|}$ \ \, and \, \ $h_m(x) \larrow \frac{1}{|x|}$, \ for a.e. $x\in \R^n$,\ \,  as $m\ri \infty$.} 
\]
The limit in \eqref{2.4} is meant in the weak $L^{2^*}$ sense as well as a.e.\ on $\R^n$, and $u$ is independent of the choice of sequence $(h_m)^\infty_1 $.

\et

\begin{remark} The solution $u$ above is vectorial but \emph{real}, although the formula \eqref{2.4} involves complex quantities. Moreover, ``$L^{2^*}$" above means ``$L^{p}$ for $p={2^*}$", not ``$(L^2)^*$".
\end{remark}

\noi \textbf{Formal derivation of the representation formula.} Before giving the rigorous proof of Theorem \ref{th1}, it is very instructive to derive \emph{formally} a representation formula for the solution of $\A  : Du =f$. By applying the Fourier transform to the PDE, we have
\[
\A : \widehat{Du}\,=\, \widehat{f}, \ \ \text{ a.e.\ on }\R^n,
\]
and hence,
\[
2\pi i\, \A : \widehat{u}(z) \ot z \,= \, \widehat{f}(z),  \ \  \text{ for a.e. }z\in \R^n.
\]
For clarity, let us also rewrite this equation in index form:
\[
\left(\A_{\al \be j} {z_j}\right)  \widehat{u_\be}(z)\,=\, \frac{1}{2\pi i}\widehat{f_\al}(z).
\]
Hence, we have
\[
\left(\A \frac{z}{|z|}\right)  \widehat{u}(z) \, = \, \frac{1}{2\pi i |z|}\widehat{f}(z)
\]
and by using the identity
\beq \label{iden}
\big(\A \sgn(z)\big)^{-1}\, =\, \frac{\ \cof \big(\A \sgn(z)\big)^\top}{ \det \big(\A \sgn(z)\big)}
\eeq
we get
\begin{align}
\widehat{u}(z) \, &=\,  \frac{ 1 }{2\pi i |z|}\big(\A \sgn(z)\big)^{-1}  \widehat{f}(z)  \nonumber\\
&=\,  \frac{1}{2\pi i |z|} \frac{\ \cof \big(\A \sgn(z)\big)^\top}{ \det \big(\A \sgn(z)\big)}\widehat{f}(z). \nonumber
\end{align}
By the Fourier inversion formula and  the identity $f^{\vee}(z)=\widehat{f}(-z)=\big(f(-\cdot)\big)^{\wedge}$, we obtain
\begin{align}
u \,&= \, \frac{1}{2\pi i} \left\{\frac{1}{|\cdot|}\frac{\ \cof \big(\A \sgn\big)^\top}{ \det \big(\A \sgn\big)} \widehat{f} \right\}^{\vee}    \nonumber\\
&= \, -\frac{1}{2\pi i } \left\{\frac{1}{|\cdot|}\frac{\ \cof \big(\A \sgn\big)^\top}{ \det \big(\A \sgn\big)}\overset{\vee}{f}\right\}^{\wedge} . \nonumber
\end{align}
Hence, we get the formula
\beq \label{2.5}
u \,= \, -\frac{1}{2\pi i} \left\{\widehat{ \frac{1}{|\cdot|}} \ast \left[ \frac{\ \cof \big(\A \sgn\big)^\top}{ \det \big(\A \sgn\big)}\overset{\vee}{f}\right]^{\wedge}\right\}.
\eeq
Formula \eqref{2.5} is ``the same" as \eqref{2.4}, \emph{if we are able to pass the limit inside the integrals} of the convolution and the Fourier transform. However, this may not be possible since $h=1/|\cdot|$ is only in the \emph{weak $L^1$} space $L^{1,\infty}(\R^n)$. Convergence needs to be rigorously justified, and this is the content of the proof of Theorem \ref{th1}. Further, by using the next identity (which follows by the properties of the Riesz potential)
\[
\left(\frac{1}{|\cdot|}\right)^{\wedge} =\, \ga_{n-1} \frac{1}{\ |\cdot|^{n-1}}
\]
where the constant $\ga_{\al}$ equals 
\[
\ga_\al \,=\, \frac{2^\al\, \pi^{n/2} \,\Gamma(\al/2)}{\Gamma(n/2-\al/2)}, \ \ \  \, 0\, <\, \al\, <\, n,
\]
we may rewrite \eqref{2.5} as
\beq \label{2.6}
u \, =\, -\frac{\ga_{n-1}}{2\pi i |\cdot|^{n-1}} \ast \left[ \frac{\ \cof \big(\A \sgn\big)^\top}{ \det \big(\A \sgn\big)}\overset{\vee}{f}\right]^{\wedge}.
\eeq
Formula \eqref{2.6} is the formal interpretation of the expression \eqref{2.4}, which we will now establish rigorously.

\BPT \ref{th1}. We begin by assuming that a solution of $\A :Du =f$ exists, and we derive the a priori estimate. By applying the Fourier transform and arguing as above, we have
\[
2\pi i\, \A : \widehat{u}(z) \ot z \, =\,  \widehat{f}(z),
\]
for a.e.\ $z\in \R^n$. Let $\nu\equiv \nu(\A)$ be the ellipticity constant of $\A$ (see \eqref{2.1a}). We then get
\[
2\pi \nu\, |\widehat{u}(z)| \, |z| \ \leq \ 2\pi \, |\A : \widehat{u}(z) \ot z| \, =\,  \big|\widehat{f}(z) \big|
\]
and hence we have
\begin{align}
\big| \widehat{Du}(z)\big|^2 \ &=\ \big|2\pi i \, \widehat{u}(z) \ot z \big|^2 \nonumber \\
& =\  \big| \widehat{u}(z)\big|^2|2\pi z|^2  \nonumber\\
& \leq \ \frac{1}{\nu^2} \big|\widehat{f}(z) \big|^2, \nonumber
\end{align}
for a.e.\ $z\in \R^n$. By integrating the above inequality on the whole space, Plancherel's theorem gives 
\[
\| Du\|_{L^2(\R^n)} \, \leq \, \frac{1}{\nu} \| f\|_{L^2(\R^n)}.
\]
Further, since $n\geq 3$, the Gagliardo-Nirenberg-Sobolev inequality gives that there exists $C=C(n,N)>0$ such that
\beq \label{2.7}
\| u\|_{L^{2^*}(\R^n)} \, \leq \, C\| Du\|_{L^2(\R^n)} \, \leq\, \frac{C}{\nu} \| f\|_{L^2(\R^n)}.
\eeq
Hence, \eqref{2.7} implies \eqref{2.3}. Now we prove existence of $u$ and the desired formula \eqref{2.5}. Let $(h_m)^\infty_1 \sub \S(\R^n)$ be any sequence of even functions in the Schwartz class for which 
\beq \label{2.7a}
\text{$0\leq h_m(x) \leq \frac{1}{|x|}$ \  and \ $h_m(x) \larrow \frac{1}{|x|}$, \ for a.e.\ $x\in \R^n$,\  as $m\ri \infty$.} 
\eeq
We set:
\beq \label{2.8}
u_m\, :=\, -\frac{1}{2\pi i} 
 \widehat{h_m} \ast
\left[ 
\frac{ \ \cof \,(\A \sgn )^\top }{ \det (\A \sgn) }\, \overset{\vee}{f}
\right]^{\wedge} .
\eeq
We will now show that the function $ u_m$ of \eqref{2.8} satisfies
\[
u_m\, \in\, L^2(\R^n)^N\cap L^\infty(\R^n)^N. 
\]
Indeed, observe first that since ${h_m} \in \S(\R^n)$ and the Fourier transform is bijective on the Schwartz class, we have
\[
\widehat{h_m}\, \in \,\S(\R^n) \, \sub \, L^1(\R^n)\cap L^2(\R^n).
\]
Let now $p\in [1,2]$ and define $r$ by
\[
r\ :=\ \frac{2p}{2-p}.
\]
Then, we have
\[
1\, +\, \frac{1}{r}\, =\, \frac{1}{p}\, +\, \frac{1}{2}, \ \ \ 1\leq p\leq2,
\]
and by Young's inequality and Plancherel's theorem, we obtain
\begin{align}
\|u_m\|_{L^r(\R^n)}\ &\leq \  \frac{1}{2\pi} \big\| \widehat{h_m} \big\|_{L^p(\R^n)} 
\left\| \left[ 
\frac{ \ \cof \,(\A \sgn )^\top }{ \det (\A \sgn) }\, \overset{\vee}{f}
\right]^{\wedge}   
\right\|_{L^2(\R^n)}   \nonumber\\
 &\leq \  \frac{1}{2\pi} \big\| \widehat{h_m} \big\|_{L^p(\R^n)} 
\left\|
\frac{ \ \cof \,(\A \sgn )^\top }{ \det (\A \sgn) }\, \overset{\vee}{f}
\right\|_{L^2(\R^n)} \nonumber .
\end{align}
We now recall that the estimate \eqref{2.2}  implies 
\[
\underset{z\in \R^n}{\ess\,\inf}\, \big|\det (\A \sgn(z))\big|\, >\, 0 
\]
and hence we get
\begin{align}
\|u_m\|_{L^r(\R^n)}\ 
&\leq \ \frac{1}{2\pi} \big\| \widehat{h_m} \big\|_{L^p(\R^n)} 
\left\|
\frac{ \, \cof \,(\A \sgn ) }{ \det (\A \sgn) }\right\|_{L^\infty(\R^n)}  
\big\|\overset{\vee}{f}
\big\|_{L^2(\R^n)}   \nonumber\\
&\leq \  C \big\| \widehat{h_m} \big\|_{L^p(\R^n)} 
\left\|f
\right\|_{L^2(\R^n)}, \nonumber
\end{align}
for some $C>0$ depending only on $|\A|$ and $\nu(\A)$. Consequently, $u_m  \in L^r(\R^n)^N$ for all $r \in [2, \infty]$.  

Next, by \eqref{2.8} and the properties of convolution, we obtain
\[
u_m\, =\, -\frac{1}{2\pi i} 
\left[ 
h_m \frac{ \ \cof \,(\A \sgn )^\top }{ \det (\A \sgn) }\, \overset{\vee}{f}
\right]^{\wedge} ,
\]
a.e.\ on $\R^n$. The Fourier inversion theorem gives
\[
\overset{\vee}{u_m}\, =\, -\frac{1}{2\pi i}  h_m 
\frac{ \ \cof \,(\A \sgn )^\top }{ \det (\A \sgn) }\, \overset{\vee}{f},
\]
a.e.\ on $\R^n$. Since $h_m(-z)=h_m(z)$ for all $z\in \R^n$, we get
\begin{align}
\widehat{u_m} (z)\, &=\, -\frac{1}{2\pi i}  h_m(z)
\frac{ \ \cof\, \Big( -\A \dfrac{z}{|z|} \Big)^\top }{ \det \Big(-\A \dfrac{z}{|z|} \Big) }\, \widehat{f}(z) \nonumber\\
&=\, -\frac{1}{2\pi i}  h_m(z)\frac{ (-1)^{N-1}}{ (-1)^{N}}
\frac{\ \cof\, \Big( \A \dfrac{z}{|z|} \Big)^\top }{ \det \Big(\A \dfrac{z}{|z|} \Big) }\, \widehat{f}(z) . \nonumber
\end{align}
Hence, by the identity \eqref{iden}, we deduce
\[
\widehat{u_m} (z)\, =\, \frac{1}{2\pi i}  h_m(z) \Big( \A \dfrac{z}{|z|} \Big)^{-1}\, \widehat{f}(z) ,
\]
a.e.\ on $\R^n$, which we rewrite as
\beq \label{2.10}
 \A :\widehat{u_m} (z) \ot 2\pi i z\, =\, \big(h_m(z)|z|\big)\, \widehat{f}(z).
\eeq
Equivalently,
\beq \label{2.10a}
\A : \widehat{Du_m}(z)\, =\, \big(h_m(z)|z|\big)\, \widehat{f}(z).
\eeq
By \eqref{2.10} we have that 
\beq \label{2.14}
0 \, \leq \, h_m(z)|z| \, \leq\, 1
\eeq
and hence by \eqref{2.14}, \eqref{2.10}, \eqref{2.10a}, \eqref{2.7a} and in view of \eqref{2.1}, we may argue again as in the derivation of \eqref{2.7} to obtain that each $u_m$ satisfies the estimate \eqref{2.3}. Hence, there is a subsequence of $m$'s and a map $u\in W^{1;2^*\!,2}(\R^n)^N$ such that, along the subsequence,
\begin{align}
 &\ u_m\lharpoonup u,\ \ \ \ \text{ in $L^{2^*}(\R^n)^N$ as } m\ri \infty, \nonumber\\
 &\ u_m\larrow u,\ \ \ \ \text{a.e.\ and in $L^{2}_{\text{loc}}(\R^n)^N$ as } m\ri \infty, \nonumber\\
 &Du_m\lharpoonup Du,\, \text{ in $L^{2}(\R^n)^{Nn}$ as } m\ri \infty. \nonumber
\end{align}
By \eqref{2.14} and since $h_m(z)|z|\ri 1$ for a.e.\ $z\in \R^n$, the Dominated Convergence theorem implies
\[
\big|h_m\, | \cdot | \big|\widehat{f} \larrow \widehat{f}, \ \text{ in  $L^2(\R^n)^N$ as } m\ri \infty. 
\]
By passing to the limit as $m\ri \infty$ in \eqref{2.10a}, since both $Du$ and $f$ are $L^2$ maps, the Fourier inversion formula implies that $u$ solves 
\[
\A:Du\, =\, f
\]
a.e.\ on $\R^n$. By passing to the limit as $m\ri \infty$ in \eqref{2.8}, we obtain the desired representation formula \eqref{2.4}. Uniqueness of the limit $u$ (and hence independence from the choice of sequence $h_m$) follows from the a priori estimate \eqref{2.3} and linearity. The theorem ensues.       \qed

\section{Strict ellipticity and Existence-uniqueness in the fully nonlinear case} \label{section3}
\ms

In this section we focus on the derivation of the appropriate condition allowing to prove existence and uniqueness of solution in the fully nonlinear case of the PDE system
\beq  \label{3.1}
F(\cdot,Du ) \,=\,  f, \ \ \ \text{ a.e.\ on }\R^n.
\eeq
Here and subsequently $F\, :\, \R^{n} \by \R^{N \by n} \larrow  \R^N$  is a Carath\'eodory map, namely
\beq \label{3.2}
\left\{
\begin{array}{l} 
x\mapsto F(x,P) \text{ is measurable, for every } P\in \R^{N\by n},\ms\\
P\mapsto F(x,P) \text{ is continuous, for a.e. } x\in \R^{n}.
\end{array}
\right.
\eeq
The crucial assumption in order to prove unique solvability of \eqref{3.1} is the next strict ellipticity condition.

\begin{definition}[Strict ellipticity] \label{def1} Let $F\, :\, \R^{n} \by \R^{N \by n} \larrow  \R^N$ satisfy \eqref{3.2}. We say  that \eqref{3.1} is an elliptic system (or that $F$ is elliptic) when there exists a linear map
\[
\A \ : \ \R^{N\by n} \larrow \R^N
\]
such that
\beq \label{3.3}
\underset{x\in \R^n}{\ess\,\sup} \sup_{P,Q\neq 0}\,\left| \frac{F(x,P+Q) -F(x,P) -\A:Q}{|Q|}\right| \ <\, \min_{|\eta|=|a|=1}\big| \A:\eta \ot a \big|.
\eeq
\end{definition}

\ms
\noi We recall that for the right hand side we have the notation $\nu(\A)$ of \eqref{2.1a}.

\begin{remark} In the sequel we will assume that 
\[
\nu(\A)\, >\, 0, 
\]
which means that the linear map $\A : \R^{N\by n}\larrow \R^N$ assumed above is elliptic in the sense of \eqref{2.1}. Otherwise, if $\nu(\A)=0$, it easy to see that we have $F(x,P)=\A:P$ and then we reduce to the linear case studied in Section \ref{section2}. 
\end{remark}

\begin{remark} Loosely speaking, the meaning of \eqref{3.3} is that the difference quotient of $F(x,\cdot)$ is uniformly close to an elliptic constant tensor $\A$, and ``how close" is determined by ``how much elliptic" $\A$ is. That is, the larger the value of the ellipticity constant $\nu(\A)$ of $\A$, the larger the deviation of $F$ from this $\A$ is allowed to be. 

In particular, in the \emph{linear non-constant case} of 
\[
F(x,P)\, =\, \A(x):P, \ \ \ \ \A \, :\, \R^n \larrow \R^N \ot \R^{N\by n}, \, \text{ measurable},
\]
which corresponds to the linear system
\beq \label{3.6}
\A(x): Du(x)\, =\, f(x),
\eeq
the ellipticity assumption \eqref{3.3} simplifies to
\beq \label{3.5}
\underset{x\in \R^n}{\ess\,\sup} \sup_{|Q|=1}\,\left| \big(\A(x)-\A \big):Q\right| \ <\, \min_{|\eta|=|a|=1}\big| \A:\eta \ot a \big|.
\eeq
Hence, by using the norm 
\[
\|\A\|\, :=\, \sup_{|Q|=1}|\A:Q| \, =\, \sup_{|\xi|=|Q|=1}\big|\A_{\al \be j} \xi_\al  Q_{\be j} \big| 
\]
on $\R^N \ot \R^{N\by n}$, assumption  \eqref{3.3} says
\beq \label{3.7}
\underset{x\in \R^n}{\ess\,\sup}\, \big\| \A(x)-\A \big\| \ <\, \min_{|\eta|=|a|=1}\big| \A:\eta \ot a \big|.
\eeq
Hence, \emph{the linear system \eqref{3.6} is elliptic when there is a constant elliptic tensor $\A$ such that the distance $\|\A(x)-\A\|$ is slightly smaller than the ellipticity constant of the tensor $\A$.}
\end{remark}

\begin{remark}
Nontrivial fully nonlinear examples of maps $F$ which are elliptic in the sense of the Definition \ref{def1} above are easy to find. Consider any fixed tensor $\A \in \R^N \ot \R^{N\by n}$ for which $\nu(\A)>0$ and any Carath\'eodory map
\[
f\ : \ \R^n \by \R^{N\by n}\larrow \R^N
\]
which is Lipschitz with respect to the second variable and whose Lipschitz constant is essentially uniformly strictly smaller than the ellipticity constant of $\A$:
\[
\big\|f(x,\cdot)\big\|_{C^{0,1}(\R^{N\by n})} \ \leq\ \la \, \nu(\A), \ \ \text{for a.e. }x\in \R^n, \ \ 0<\la<1.
\]
Then, the map $F :\R^n \by \R^{N\by n}\larrow \R^N$ given by
\[
F(x, Q)\, :=\, \A:Q\, +\, f(x, Q)
\]
satisfies
\begin{align}
\Big|F(x,P+Q)-F(x,P)-\A:Q \Big|\, &=\, \big|f(x,P+Q)-f(x,P) \big| \nonumber\\
&\leq\, \la \,\nu(\A) |Q|,\nonumber
\end{align}
and hence is elliptic in the sense of \eqref{3.3}. \ms

\noi \textit{Thus, every Lipschitz perturbation of an elliptic constant tensor gives a fully nonlinear elliptic map, when the Lipschitz constant of the perturbation is strictly smaller than the ellipticity constant of the tensor.}

\end{remark}

We now show that the ellipticity assumption can be seen an a notion of pseudo-monotonicity, coupled by Lipschitz continuity of $Q\mapsto F(x,Q)$. 

\begin{lemma}[Relation of ellipticity and pseudo-monotonicity] \label{pr1} Suppose that the map $F : \R^n \by \R^{N\by n}\larrow \R^N$ satisfies \eqref{3.2}.
Consider the statements
\begin{enumerate}

\item There exists $\A \in \R^N \ot \R^{N\by n}$ with $\nu(\A)>0$ such that $F$ is strictly elliptic, namely satisfies the inequality \eqref{3.3}.

\ms

\item \begin{itemize}
         \item $Q\mapsto F(x,Q)$ is globally Lipschitz continuous on $\R^{N\by n}$, essentially uniformly in $x\in \R^n$.
         
         \item \emph{(Pseudo-Monotonicity)} There exists an $\A \in \R^N \ot \R^{N\by n}$ for which $\nu(\A)>0$ and also a $\la\in(0,1)$ such that, for all $P,Q \in \R^{N \by n}$ and a.e.\ $x\in \R^n$, 
          \begin{align} \label{3.8}
       \ \ \ \  \ \ \ \ \  (\A:Q)^\top\Big[F(x,P+Q) - F(x,P) \Big]\ \geq  \ \frac{1}{2}|\A:Q|^2\, -\, \frac{\la^2}{2}\nu(\A)^2|Q|^2,
            \end{align}
where $\nu(\A)$ is given by \eqref{2.1a}.
        \end{itemize}

\end{enumerate}
\ms

Then, $(1)$ implies $(2)$. Conversely, $(2)$ implies $(1)$ when in addition the Lipschitz constant of $Q\mapsto F(x,Q)$ is small enough:
\beq  \label{3.8a}
\underset{x\in \R^n}{\ess\, \sup}\, \big\|F(x,\cdot)\|_{C^{0,1}(\R^{N\by n})}\ <\ \sqrt{1-\la^2}\, \nu(\A).
\eeq
\end{lemma}

\BPL \ref{pr1}. Assume $(1)$. By \eqref{3.3} we have
\[
\big|F(x,P+Q)-F(x,P) \big|\, \leq\, \big(\nu(\A) +\|\A\|\big)|Q|,
\]
for a.e.\ $x\in \R^n$ and all $P,Q\in \R^{N\by n}$. Hence, $F(x,\cdot)$ is Lipschitz, essentially uniformly in $x$. Again by \eqref{3.3}, we have that there exists $\la \in (0,1)$ such that
\[
\Big| F(x,P+Q) -F(x,P) -\A:Q \Big| \leq \, \la \nu(\A)|Q|.
\]
Hence, 
\begin{align}  \nonumber
 \la^2 \nu(\A)^2|Q|^2 \, &\geq \,  \big|F(x,P+Q) -F(x,P) \big|^2+\, | \A: Q |^2  \nonumber \\
& \ \ \ -\, 2\left( \A: Q \right)^\top\Big[  F(x,P+Q) -F(x,P) \Big] \nonumber\\
&\geq\, | \A: Q |^2-\, 2\left( \A: Q \right)^\top\Big[  F(x,P+Q) -F(x,P) \Big]. \nonumber
\end{align}
The above inequality implies \eqref{3.8}, and $(2)$ ensues.  Conversely, assume $(2)$ and also \eqref{3.8a}. Then, by \eqref{3.8a} there exists $\de \in (0,1)$ such that
\[
\big|F(x,P+Q) -F(x,P) \big|^2 \leq\,  \de^2 \big( 1-\la^2 \big)\nu(\A)^2|Q|^2.
\]
By adding this inequality to 
\[
|\A:Q|^2\,-\, 2 (\A:Q)^\top\Big[F(x,P+Q) - F(x,P) \Big]\ \leq  \, \la^2\nu(\A)^2|Q|^2
\]
we get
\[
\Big| F(x,P+Q) -F(x,P) -\A:Q \Big| \leq \, \sqrt{\la^2 +\, \de^2 \big(1-\la^2\big)}\, \nu(\A)|Q|.
\]
Since $\la^2 +\, \de^2 (1-\la^2)<1$, we see that the above inequality implies \eqref{3.3} and hence $(1)$ ensues, as desired.    \qed

\ms

The main result of this paper is the next theorem:

\bt[Existence-Uniqueness] \label{th2} Assume that $n\geq 3$, $N\geq 2$ and let $F : \R^n \by \R^{N \by n}\larrow \R^N$ be a Carath\'eodory map, satisfying \eqref{3.3} and also $F(x,0)=0$ for a.e.\ $x\in \R^n$. Let also $f\in L^2(\R^n)^N$. Then, the problem
\[
F(\cdot,Du)\, =\, f, \ \ \text{ a.e.\ on }\R^n,
\]
has a unique solution $u$ in the space $W^{1;2^*\!,2}(\R^n)^N$ (see \eqref{1.5}), which also satisfies the estimate
\beq \label{4.1}
\|u\|_{L^{2^*}(\R^n)}\ + \ \|Du\|_{L^{2}(\R^n)}\, \leq\, C \|f\|_{L^{2}(\R^n)}
\eeq
for some $C>0$ depending only on $F$. 
\et
 
In the course of the proof we will establish the following strong uniqueness estimate, which is a form of ``comparison principle in integral norms":

\bcor[Uniqueness estimate] \label{cor1} Assume that $n\geq 3$, $N\geq 2$ and let $F : \R^n \by \R^{N \by n}\larrow \R^N$ be a Carath\'eodory map, satisfying \eqref{3.3}. Then, for any two maps $w,v \in W^{1;2^*\!,2}(\R^n)^N$, we have
\beq \label{4.2}
\|w-v\|_{L^{2^*}(\R^n)}\ + \ \|Dw-Dv\|_{L^{2}(\R^n)}\, \leq\, C \big\|F(\cdot,Dw)-F(\cdot,Dv) \big\|_{L^{2}(\R^n)}.
\eeq
In particular, any two global strong a.e.\ solutions of the PDE system $F(\cdot,Du)\, =\, f$ coincide.
\ecor

The proofs of Theorem \ref{th2} and Corollary \ref{cor1} utilise the following result of Campanato taken from \cite{C0}, whose short proof is given for the sake of completeness at the end of the section:

\bt[Campanato's near operators] \label{th3}  Let $F,A : \mathfrak{X} \larrow X$ be two maps from the set $\mathfrak{X} \neq \emptyset$ to the Banach space $(X,\|\cdot\|)$. Suppose there exists $0<K<1$ such that
\beq \label{4.3}
\Big\|F[u]-F[v]-\big( A[u]-A[v]\big) \Big\| \, \leq\, K \big\| A[u]-A[v] \big\|,
\eeq
for all $u,v \in \mathfrak{X}$. Then, if $A$ is a bijection, $F$ is a bijection as well.
\et
Campanato defined the inequality \eqref{4.3}  above as the ``nearness of $F$ to $A$", using also a multiplicative constant of front of (either $A$ or) $F$. Such a constant has no bearing in the generality we are working in, so we normalise it to one in the definition.

\BPT \ref{th2} (\textbf{and Corollary} \ref{cor1}). By our assumption \eqref{3.3} on $F$ and that $F(x,0)=0$, Lemma \ref{pr1} implies that there exists an $M>0$ depending only on $F$, such that for any $u\in W^{1;2^*\!,2}(\R^n)^N$, we have
\begin{align} \label{4.5}
\big\|F(\cdot,Du) \big\|_{L^2(\R^n)} \, &\leq\, \big\| F(\cdot,0)\big\|_{L^2(\R^n)}\, +\, M \| Du\|_{L^2(\R^n)} \\
&=\, M \| Du\|_{L^2(\R^n)} \nonumber\\
&\leq\, M \Big( \| Du\|_{L^2(\R^n)} \, +\, \| u\|_{L^{2^*}(\R^n)} \Big). \nonumber
\end{align}
Let also $\A \in \R^N \ot \R^{N\by n}$ be the tensor given by assumption \eqref{3.3}, which satisfies $\nu(\A)>0$, with $\nu(\A)$ as in \eqref{2.1a}. Then, we have
\begin{align} \label{4.6}
\|\A:Du\|_{L^2(\R^n)}\, &\leq\, \|\A\|\, \| Du\|_{L^2(\R^n)} \nonumber\\
&\leq\,  \|\A\| \Big( \| Du\|_{L^2(\R^n)} \, +\, \| u\|_{L^{2^*}(\R^n)} \Big).\nonumber
\end{align}
By \eqref{4.5} and \eqref{4.6} we obtain that the operators 
\[
\left\{
\begin{array}{l}
A[u]\ :=\ \A :Du, \ms\\
F[u]\ :=\ F(\cdot, Du),
\end{array}
\right.
\]
map $W^{1;2^*\!,2}(\R^n)^N$ into  $L^2(\R^n)^N$. If $u,v\in W^{1;2^*\!,2}(\R^n)^N$, then \eqref{2.1} and Plancherel's theorem give (below we denote the identity map by ``$Id$", that is $Id(x):=x$):
\begin{align} \label{4.6}
\big\|\A:Du -\A:Dv \big\|_{L^2(\R^n)}\, &= \, \big\|\A:\widehat{Du} - \A:\widehat{Dv} \big\|_{L^2(\R^n)} \nonumber\\
&= \, \big\| \A: \big(\widehat{u}-\widehat{v}\big) \ot (2\pi i Id) \big\|_{L^2(\R^n)}   \nonumber \\
&\geq \, \nu(\A)\big\| \big(\widehat{u}-\widehat{v}\big) \ot (2\pi i Id) \big\|_{L^2(\R^n)} \\
&= \, \nu(\A)\big\|\widehat{Du} - \widehat{Dv} \big\|_{L^2(\R^n)}  \nonumber\\
&= \, \nu(\A)\big\| {Du} -  {Dv} \big\|_{L^2(\R^n)}.   \nonumber
\end{align}
We now set
\[
\nu(F,\A)\, :=\, \underset{x\in \R^n}{\ess\,\sup} \sup_{P,Q\neq 0}\,\left| \frac{F(x,P+Q) -F(x,P) -\A:Q}{|Q|}\right| .
\]
In view of \eqref{2.1a}, we may rewrite  \eqref{3.3}  as
\beq \label{3.15a}
0\,<\,  \nu(F,\A)\, <\, \nu(\A).
\eeq
By employing \eqref{3.3}, for $u,v \in W^{1;2^*\!,2}(\R^n)^N$ we have
\begin{align} 
\Big\|F(\cdot,Du) &-  F(\cdot,Dv)  -\A: \big(Du-Dv \big)\Big\|_{L^2(\R^n)} \nonumber\\
&\ \leq \left(\underset{\R^n}{\ess\,\sup} \sup_{P,Q\neq 0}\,\left| \frac{F(\cdot,P+Q) -F(\cdot,P) -\A:Q}{|Q|}\right| \right)
\big\| Du-Dv\big\|_{L^2(\R^n)}  \nonumber\\
&\ = \, \nu(F,\A) \big\| Du-Dv\big\|_{L^2(\R^n)} \nonumber\\
&\overset{\eqref{4.6}}{\leq}\, \frac{\nu(F,\A)}{\nu(\A)} \big\|\A:(Du -Dv)\big\|_{L^2(\R^n)} \nonumber
\end{align}
and hence we obtain the inequality
\beq   \label{4.7}
\Big\|F(\cdot,Du) -  F(\cdot,Dv)  -\A: \big(Du-Dv \big)\Big\|_{L^2(\R^n)}  \leq \, \frac{\nu(F,\A)}{\nu(\A)} \big\|\A:(Du -Dv)\big\|_{L^2(\R^n)}  .
\eeq
We now recall that since $\nu(\A)>0$, Theorem \ref{th1} implies that the linear operator
\[
A\ :\  W^{1;2^*\!,2}(\R^n)^N \larrow  L^2(\R^n)^N 
\]
is a bijection. Hence, in view of the inequalities \eqref{3.15a}  and \eqref{4.7}, Campanato's Theorem \ref{th3} implies that $F$ is a bijection as well. As a result, for any $f\in L^2(\R^n)^N$, the PDE system 
\[
F(\cdot, Du)\, =\, f,\ \  \text{ a.e. on }\R^n,
\]
has a unique solution $u\in W^{1;2^*\!,2}(\R^n)^N$. Moreover, by \eqref{4.7} we deduce the estimate
\begin{align} 
\Big\|F(\cdot,Du) -  F(\cdot,Dv)\Big\|_{L^2(\R^n)} \,  &\geq \, \left(1-\frac{\nu(F,\A)}{\nu(\A)} \right) \big\|\A:(Du -Dv)\big\|_{L^2(\R^n)}\nonumber \\
&\geq\, \big(\nu(\A)-\nu(F,\A) \big)\, \big\| Du -Dv\big\|_{L^2(\R^n)} .\nonumber
\end{align}
This last estimate together with the fact that $n\geq 3$ and the Gagliardo-Nirenberg-Sobolev inequality, imply both \eqref{4.1} and  \eqref{4.2}. The theorem ensues, and so does Corollary \ref{cor1}.      \qed

\ms

We conclude this section with the proof of Campanato's theorem on near operators taken from \cite{C0}, which we provide for the convenience of the reader. 

\BPT \ref{th3}. It suffices to show that for any $f\in X$, there is a unique $u\in \mathfrak{X}$ such that
\[
F[u]\, =\, f.
\]
In order to prove that, we first turn $\mathfrak{X}$ into a complete metric space, by pulling back the structure from $X$ via $A$: for, we define the distance
\[
d(u,v)\, :=\, \big\| A[u]-A[v]\big\|.
\]
Next, we fix an $f\in X$ and define the map
\[
T\ : \ \mathfrak{X} \larrow \mathfrak{X}\ , \ \ \ T[u]\, :=\, A^{-1}\Big(A[u]-\big(F[u]-f \big) \Big).
\]
We conclude by showing that $T$ is a contraction on $(\mathfrak{X},d)$, and hence has a unique $u\in \mathfrak{X}$ such that $T[u]=u$. The latter equality is equivalent to $F[u]=f$, and then we will be done.  Indeed, we have that
\begin{align}
d\Big( T[u],T[v] \Big)  \, &=\ \Big\|  \left(A[u]-\big(F[u]-f \big) \right) \, -\,  \left(A[v]-\big(F[v]-f \big) \right) \Big\| \nonumber \\
 &=\ \Big\| A[u] -A[v] -\big(F[u]  - F[v]\big) \Big\|, \nonumber 
\end{align}
and hence
\begin{align}
d\Big( T[u],T[v] \Big)  \ \, 
&\!\!\!\!\overset{\eqref{4.3}}{\leq}   K \big\| A[u]-A[u] \big\| \nonumber\\
  &= \  K \, d(u,v). \nonumber
\end{align}
Since $K<1$, the conclusion follows and the theorem ensues.              \qed

\ms

\ms

\noi \textbf{Acknowledgement.} I would like to thank Jan Krinstensen and Bernard Dacorogna for their share of expertise on the status of the subject. I am also indebted to Tristan Pryer and Beatrice Pelloni for our inspiring scientific discussions. Finally, I would like to thank the referee of this paper who's suggestions improved the content and the appearance of this work.

\ms

\ms

\bibliographystyle{amsplain}

\end{document}